\documentclass[a4paper,12pt]{article}
\usepackage{latexsym}
\usepackage{amsmath}
\usepackage{amssymb}
\usepackage{enumerate}

\usepackage{theorem}
\newtheorem{theorem}{Theorem}[section]

\newtheorem{lemma}[theorem]{Lemma}
\newtheorem{corollary}[theorem]{Corollary}

\newtheorem{MT}{Main Theorem.}

\newtheorem{CMT}{Corollary to Main Theorem.}

\theorembodyfont{\rmfamily}
\newtheorem{proof}{\textmd{\textit{Proof.}}}

\newtheorem{remark}[theorem]{Remark}
\newtheorem{example}[theorem]{Example}

\newtheorem{acknowledgement}{\textmd{\textit{Acknowledgement.}}}

\makeatletter

\@addtoreset{equation}{section}
\makeatother

\newcommand{\qedd}{\hfill \Box}
\newcommand{\ve}{\varepsilon}
\newcommand{\lra}{\longrightarrow}
\newcommand{\wt}{\widetilde}
\newcommand{\wh}{\widehat}

\newcommand{\R}{\ensuremath{\mathbb{R}}}
\newcommand{\Sph}{\ensuremath{\mathbb{S}}}

\newcommand{\cE}{\ensuremath{\mathcal{E}}}

\newcommand{\cN}{\ensuremath{\mathcal{N}}}

\newcommand{\cU}{\ensuremath{\mathcal{U}}}

\def\Cut{\mathop{\mathrm{Cut}}\nolimits}

\title
{
Total Curvatures of Model Surfaces Control\\ 
Topology\hspace{-0.5mm} of Complete Open Manifolds with\\ 
Radial Curvature Bounded Below.\,I\footnote{
Mathematics Subject Classification (2010)\,:\,53C21, 53C22.}
\footnote{
Keywords\,: Busemann function, radial curvature, total curvature of surface}
}
\author{Kei KONDO $\cdot$ Minoru TANAKA}
\date{}
\begin{document}
\maketitle

\begin{abstract}
We investigate the finiteness structure of 
a complete non-compact $n$-dimensional Riemannian manifold $M$ whose radial curvature at 
a base point of $M$ is bounded from below by that of a non-compact von Mangoldt surface of 
revolution with its total curvature greater than $\pi$.
We show, as our main theorem, that 
all Busemann functions on $M$ are exhaustions, and that 
there exists a compact subset of $M$ such that the compact set contains 
all critical points for any Busemann function on $M$. 
As corollaries by the main theorem, 
$M$ has finite topological type, and the isometry group of $M$ is compact. 
\end{abstract}
\section{Introduction}\label{sec:int}

\noindent

There is great interest in the relationship between 
radial (sectional) curvature geometry and pure sectional curvature geometry. 
The difficulty in the geometry of radial curvatures from a given point  is that 
one has to look around such manifolds only at the base point. 
For example, in the Toponogov comparison theorem, 
(See Theorem \ref{TCT} in this article), in such a comparison geometry, 
all geodesic triangles must have the base point as one of their vertices. 
Thus, radial curvature geometry does not have, in a sense, 
{\em homogeneity on curvature} like pure sectional curvature geometry, 
which has all curvatures everywhere bounded below by some constant 
so that the Toponogov comparison theorem in such a pure geometry 
holds for all geodesic triangles. 
Despite such a difficulty in radial curvature geometry, 
we have seen the new results (See \cite[Gap Theorem]{GW}, \cite{KK}). 
In particular, after the work \cite{IMS}, and also \cite{SiT}, 
the difficulty in radial curvature geometry increases more, 
since von Mangoldt surfaces of revolution (defined below), 
and $2$-spheres of revolution, and more general classes of models are employed 
as reference spaces in comparison theorems of radial curvature 
geometry (Note that Hadamard surfaces with finite total curvature have been 
employed as reference models in \cite{GW}, \cite{A}, and other articles). 
However, we have obtained some results as the relationships between 
radial curvature geometry and pure sectional curvature geometry in 
\cite{KO}, \cite{K}, \cite{ST}, \cite{KT1}, \cite{KT2}, and \cite{KT3}.

\bigskip

The Gauss\,--\,Bonnet Theorem says that 
the total curvature of a compact Riemannian $2$-manifold $S$ 
is a topological invariant, that is, equal to $2 \pi \chi (S)$. 
Here, $\chi (S)$ is the Euler characteristic of $S$. 
For a complete non-compact Riemannian $2$-manifold, 
however, admitting a total curvature, 
the total curvature of the surface is not a topological invariant anymore. 
Cohn\,-Vossen proved that 
if a complete non-comact, finitely-connected Riemannian $2$-manifold $X$ admits a total curvature, 
then the total curvature of $X$ is not greater than $2 \pi \chi (X)$ (See \cite[Satz 6]{CV1}). 
He has developed fundamental techniques for investigating 
the structures of complete non-compact Riemannian $2$-manifolds. 
Although he restricted himself to $2$-dimensional Riemannian manifolds, 
some techniques are even now useful for investigating the relationship between 
the topology and the sectional curvature of higher-dimensional complete non-compact 
Riemannian manifolds. As pointed out in the preface of \cite{SST}, 
{\em 
it took more than thirty years to obtain higher-dimensional extensions of 
Cohn\,-Vossen's results.} 
They are the Toponogov splitting theorem \cite{To} and the structure theorem for 
complete non-compact Riemannian $n$-manifolds of positive sectional curvature \cite{GM} and of 
non-negative sectional curvature \cite{CG2}. 
As pointed out above, $2 \pi \chi (X) - c(X)$, where $c(X)$ denotes the total curvature of 
a complete non-compact, finitely-connected Riemannian $2$-manifold $X$, 
is not a topological invariant, but it depends only on the ends of $X$. 
This is a direct consequence of the isoperimetric inequalities (See \cite[Theorem 5.2.1]{SST}).

\bigskip

Our main purpose in this article, from radial curvature geometry's standpoint, 
is to generalize the following result of Shiohama in the geometry of complete non-compact surfaces 
to $n$-dimensional complete non-compact Riemannian manifolds\,:

\begin{theorem}{\rm (\cite[Main Theorem]{S1})}\label{S}
Let $X$ be a connected, complete non-comact, finitely-connected and 
oriented Riemannian $2$-manifold with one end. 
If the total curvature of $X$ is greater than $(2 \chi (X) - 1) \pi$, 
then all Busemann functions on $X$ are exhaustions. 
In particular, if the total curvature of $X$ is greater than $\pi$, 
then $X$ is homeomorphic to $\R^2$ 
and also all Busemann functions are exhaustions. 
Here, a function $F : X \lra \R$ is called an exhaustion, 
if $F^{-1}( - \infty, a ]$ is compact for all $a \in \R$. 
\end{theorem} 
In Section \ref{dbf} of this article, one can find 
the definition and some properties of a Busemann function 
on an arbitrary complete non-compact Riemannian manifold. 
Note that it is also proved in \cite[Main Theorem]{S1} that 
all Busemann functions on $X$ are non-exhaustions, 
if the total curvature of $X$ is less than $(2 \chi (X) - 1) \pi$.   

\bigskip

We will now introduce the radial curvature 
geometry for pointed complete non-compact Riemannian manifolds\,: 
Let $\wt{M}$ denote a complete $2$-dimensional Riemannian manifold 
homeomorphic to $\R^{2}$ with a base point $\tilde{p} \in \wt{M}$. 
Then, we call the pair $(\wt{M}, \tilde{p})$ a {\em non-compact model surface of revolution} 
if its Riemannian metric $d\tilde{s}^2$ is expressed 
in terms of geodesic polar coordinates around $\tilde{p}$ as 
\[
d\tilde{s}^2 = dt^2 + f(t)^2d \theta^2, \quad 
(t,\theta) \in (0,\infty) \times {\Sph_{\tilde{p}}^1}_. 
\]
Here $f : (0, \infty) \lra \R$ is a positive smooth function 
which is extensible to a smooth odd function around $0$, 
and 
$\Sph^{1}_{\tilde{p}} := \{ v \in T_{\tilde{p}} \wt{M} \ | \ \| v \| = 1 \}$. 
The function $G \circ \tilde{\gamma} : [0,\infty) \lra \R$ is called the 
{\em radial curvature function} of $(\wt{M}, \tilde{p})$, 
where we denote by $G$ the Gaussian curvature of $\wt{M}$, 
and by $\tilde{\gamma}$ any meridian emanating from 
$\tilde{p} = \tilde{\gamma} (0)$. 
Remark that $f$ satisfies the differential equation 
$f''(t) + G (\tilde{\gamma}(t)) f(t) = 0$ 
with initial conditions $f(0) = 0$ and $f'(0) = 1$. For each constant number $\delta >0$, 
a sector $\wt{V} (\delta) \subset \wt{M}$ is defined by   
$\wt{V} (\delta) : = \{ \tilde{x} \in \wt{M} \, | \, 0 < \theta(\tilde{x}) < \delta\}$.
Notice that 
the $n$-dimensional model surfaces of revolution are defined similarly, and they are 
completely classified in \cite{KK}. 
The total curvature $c (\wt{M})$ of $(\wt{M}, \tilde{p})$ is formally 
defined as the improper integral, i.e., 
\[
c(\wt{M}):=\int _{\wt{M}} G_{+} \circ t\, d\wt{M} +\int_{\wt{M}} G_{-} \circ t \, d\wt{M}
\]
if 
\[
\int_{\wt{M}} G_{+}\circ t\, d\wt{M}<\infty, \quad 
{\rm or} \quad \int_{\wt{M}}G_{-}\circ t\,d\wt{M}>-\infty.
\]
Here we set 
$G_{+} (t ):=\max\{ G(\tilde{\gamma} (t)), 0 \}$
and
$G_{-} (t) :=\min\{ G(\tilde{\gamma} (t)), 0 \}$. Notice that $G = G_{+}\circ t + G_{-}\circ t$. 
If $c(\wt{M})$ exists, $c(\wt{M}) = 2 \pi (1 - \lim_{t \to \infty} f'(t))$ holds, 
since $d\wt{M} = f dt d\theta$ and $f'(0) = 1$. 
Remark that $c(\wt{M}) \le 2 \pi$ holds.\par
In this article, we are going to employ a von Mangoldt surface of revolution 
as a reference space, which is, by definition, a non-compact model surface of revolution 
whose radial curvature function is non-increasing on $[0, \infty)$. 
The cut locus $\Cut(\tilde{z})$ to each point $\tilde{z} \in \wt{M} \setminus \{ \tilde{p} \}$
of a non-compact von Mangoldt surface of revolution $(\wt{M}, \tilde{p})$ is either an empty set, 
or a ray properly contained in the meridian $\theta^{-1} (\theta (\tilde{z}) + \pi)$ lying opposite to 
$\tilde{z}$, and that the endpoint of $\Cut(\tilde{z})$ is the first conjugate point to $\tilde{z}$ 
along the minimal geodesic from $\tilde{z}$ sitting in 
$ \theta^{-1} (\theta (\tilde{z})) \cup \theta^{-1} (\theta (\tilde{z}) + \pi)$ (See \cite[Main Theorem]{T}). 
Hence, {\em any non-compact von Mangoldt surface of revolution has no pair of cut points 
in the sector $\wt{V} (\pi)$}.  
Paraboloids and $2$-sheeted hyperboloids are typical examples of a von Mangoldt surface of revolution. An atypical example of a von Mangoldt surface of revolution is the following, 
where its radial curvature function changes signs on $[0, \infty)$:

\begin{example}{\bf (Sinclair)}\label{exa1.2}
We define $f (t) := e^{- t^{2}} \tanh t$ on $[0, \infty)$. It is clear that $f$ satisfies $f(0) = 0$ 
and $\displaystyle{\lim_{t \to \infty}f (t) = 0}$. 
Moreover, we have 
\begin{align}
&{f'(t) = \frac{1}{e^{t^{2}}} \left( -2t \tanh t + \frac{1}{\cosh^{2} t} \right)}_, \label{exa1.2-1}\\
&{f''(t) = \left(4t^{2} - 2 - \frac{2}{\cosh^{2}t} \right)f(t) - \frac{4t}{e^{t^{2}}\cosh^{2}t}}_. \label{exa1.2-2} 
\end{align}
Thus, by (\ref{exa1.2-1}), we see $f'(0) = 1$ and $\displaystyle{\lim_{t \to \infty}f'(t) = 0}$. 
Furthermore, by (\ref{exa1.2-2}), we get 
\[
{G (t) := -\,\frac{f''(t)}{f(t)} = \frac{8t}{\sinh 2t} + \frac{2}{\cosh^{2}t} -4 t^{2} + 2}_.
\]
Then, we see 
\[
{\frac{d}{dt}G(t) = \frac{8\sinh 2t -16t \cosh 2t}{\sinh^{2} 2t} - \frac{2\sinh 2t}{\cosh^{2}t} -8 t < 0}
\]
on $(0, \infty)$. 
Thus, $G$ is strictly monotone decreasing on $(0, \infty)$, 
and satisfies 
$\lim_{t \downarrow 0}G(t) = 8$ and $\lim_{t \to \infty}G(t) = - \infty$.
Therefore, 
a complete non-compact Riemannian $2$-manifold $(\wt{M}, \tilde{p})$ 
with a base point $\tilde{p}$ and 
$d\tilde{s}^2 = dt^2 +  f(t)^2d \theta^2$, 
$(t,\theta) \in (0,\infty) \times \Sph_{\tilde{p}}^1$, 
is a non-compact von Mangoldt surface of revolution, 
and its radial curvature function $G$ changes signs on $[0, \infty)$. 
In particular, since $f$ satisfies $\lim_{t \to \infty} f'(t) = 0$, 
the total curvature of this $(\wt{M}, \tilde{p})$ 
is equal to $2 \pi$\,(Indeed, one may prove this by calculating the total curvature, 
or by the isoperimetric inequality in \cite[(5.2.2) of Theorem 5.2.1]{SST}). 
Other examples of $(\wt{M}, \tilde{p})$ are found in \cite{T}.    
\end{example}

\medskip

Let $(M,p)$ be a complete non-compact $n$-dimensional Riemannian manifold 
with a base point $p \in M$. We say that $(M, p)$ has 
{\em 
radial curvature at the base point $p$ bounded 
from below by that of 
a non-compact model surface of revolution $(\wt{M}, \tilde{p})$} 
if, along every unit speed minimal geodesic $\gamma: [0,a) \lra M$ 
emanating from $p = \gamma (0)$, 
its sectional curvature $K_M$ satisfies
\[
K_M(\sigma_{t}) \ge G (\tilde{\gamma}(t))
\]
for all $t \in [0, a)$ and all $2$-dimensional linear spaces 
$\sigma_{t}$ spanned by $\gamma'(t)$ 
and a tangent vector to $M$ at $\gamma(t)$. 
Notice that, if the Riemannian metric of $\wt{M}$ is $dt^2 + t^{2}d \theta^2$, or $dt^2 + \sinh^{2} t\,d \theta^2$, 
then $G (\tilde{\gamma}(t)) = 0$, or $G (\tilde{\gamma}(t)) = -1$, respectively.\par 
For this definition, the radial curvature geometry looks artificial, {\bf but this is not the case}, i.e., 
we can construct a model surface of revolution for 
any complete Riemannian manifold with an arbitrary given point 
as a base point (See \cite[Lemma 5.1]{KT1}). 
The existence of a $(\wt{M}, \tilde{p})$ is therefore {\em very natural} on the above definition.

\bigskip

Our main theorem is now stated as following\,: 
 
\begin{MT} Let $(M,p)$ be a complete non-compact Riemannian $n$-manifold $M$ 
whose radial curvature at the base point $p$ is bounded from below by
that of a non-compact von Mangoldt surface of revolution $(\wt{M}, \tilde{p})$.
If $(\wt{M}, \tilde{p})$ admits $c (\wt{M}) > \pi$, 
\begin{enumerate}[{\rm ({MT--}1)}]
\item
all Busemann functions on $M$ are exhaustions, and  
\item
there exists a compact subset $C$ of $M$ such that $C$ contains 
all critical points for any Busemann function on $M$.
\end{enumerate}
\end{MT}

\medskip\noindent
Notice that, under the assumptions in the Main Theorem, $M$ has just one end 
by \cite[(C\,--\,i) in Theorem C]{KO}. 
A generalization of the (MT--1) to an $M$ which is 
not less curved than a more general model surface of revolution has been discussed in \cite{KT2}. 
Furthermore, it follows from the Main Theorem that we have the following corollary\,:

\medskip

\begin{CMT} Let $(M,p)$ be a complete non-compact Riemannian $n$-manifold $M$ 
whose radial curvature at the base point $p$ is bounded from below by
that of a non-compact von Mangoldt surface of revolution $(\wt{M}, \tilde{p})$.
If $(\wt{M}, \tilde{p})$ admits $c (\wt{M}) > \pi$,
\begin{enumerate}[{\rm ({C\,--\,}1)}]
\item
$M$ has finite topological type, that is, 
$M$ is homeomorphic to the interior of a compact manifold with boundary. 
\item
The isometry group of $M$ is compact. 
\end{enumerate}
\end{CMT}

\begin{remark}
A related result for the (C\,--\,1), 
but for a complete non-compact Riemannian $n$-manifold with 
non-negative sectional curvature everywhere, has been obtained in 
Gromov's \cite[Subsection 1.5]{Gv}. 
On the other hand, 
a related result for the (C\,--\,2), 
but for a complete non-compact Riemannian $n$-manifold with 
non-negative sectional curvature everywhere, has been obtained in \cite[Corollary 6.2]{CG2}. 
Another related result for the (C\,--\,2) is Yamaguchi's \cite[Theorem B]{Y}, 
where he has proved that if a complete non-compact Riemannian $n$-manifold 
admits a strictly convex function with compact levels and no minimum, 
then the isometry group of such a manifold is compact. 
We have generalized  the (C\,--\,1) to an $M$ which is 
not less curved than a more general model surface of revolution $\wt{M}$
admitting just finite total curvature, i.e., $c(\wt{M}) > - \infty$ (See \cite[Main Theorem]{KT1}). 
\end{remark}

\begin{remark}
If a non-compact model surface of revolution $\wt{M}$ admits a finite total curvature, 
then, for each $\ve > 0$, there exists a compact subset $\wt{K}_{\ve}$ of $\wt{M}$ such that 
\[
{
\int_{\wt{M} \setminus \wt{K}_{\ve}} |G|\, d\wt{M}<\ve
}_.
\]
Hence, we might conjecture that the Gaussian curvature of $\wt{M}$ should be almost flat outside 
of a compact subset of $\wt{M}$. The following theorem shows that this conjecture is {\bf false} 
and that the radial curvature function $G(t)$ may change signs wildly:

\begin{theorem}{\rm (\cite{KT5})}
Let $(\wt{M}, \tilde{p})$ be a non-compact model surface of revolution. 
If $\wt{M}$ admits $- \infty < c(\wt{M}) < 2 \pi$,
then, for any $\ve > 0$, 
there exists a model surface of revolution $(\wh{M}, \wh{p}\,)$ 
with its metric 
\[
\wh{g} 
= 
dt^2 + m(t)^2d \theta^2, \quad 
(t,\theta) \in (0,\infty) \times {\Sph_{\wh{p}}^1}_, 
\]
satisfying the differential equation $m''(t) + \wh{G} (t) m(t) = 0$ 
with initial conditions $m(0) = 0$ and $m'(0) = 1$, 
and admitting a finite total curvature $c(\wh{M})$ such that 
\begin{enumerate}[{\rm (1)}]
\item
$\displaystyle{\left\|\,
G(\tilde{\gamma} (t)) - \wh{G} (t)\,
\right\|_{L_{2}} \le \ve
}$, 
\item
$\displaystyle{
c(\wt{M}) \ge c(\wh{M}) \ge c(\wt{M}) - \ve
}$ 
(respectively 
$\displaystyle{
c(\wt{M}) + \ve \ge c(\wh{M}) \ge c(\wt{M})
}$),
\item
$\displaystyle{
G(\tilde{\gamma} (t)) \ge \wh{G}(t)
}$
(respectively 
$\displaystyle{
\wh{G}(t) \ge G(\tilde{\gamma} (t))
}$) on $[0, \infty)$, and 
\item
$\displaystyle{
\liminf_{t \to \infty} \wh{G} (t) = - \infty
}$
(respectively 
$\displaystyle{
\limsup_{t \to \infty} \wh{G} (t) = \infty}$).
\end{enumerate}
\end{theorem}
\end{remark}

\bigskip

In the following sections, 
all geodesics will be normalized, unless otherwise stated.  

\begin{acknowledgement}
We would like to thank Professor Robert Sinclair 
for his atypical example of a non-compact von Mangoldt surface of revolution. 
\end{acknowledgement}

\section{Differentiability of Busemann Functions}\label{dbf}

\noindent

In this section, let $M$ denote an arbitrary complete non-compact Riemannian manifold 
without curvature assumptions. 
A \!{\em Busemann function} $F_{\gamma} : M \lra \R$ of a ray $\gamma$ on $M$ is defined by 
\[
F_{\gamma}(x) := \lim_{t \to \infty} \left\{ t - d(x, \gamma(t)) \right\}_. 
\]
$F_{\gamma}$ was first introduced by H.\,Busemann\,(\cite[Section 22]{B}) to investigate 
parallels for straight lines on a straight $G$-space where every two points are joined by a unique 
geodesic realizing the distance. 
By the definition of $F_{\gamma}$, 
we have that $|F_{\gamma}(x) - F_{\gamma}(y)| \le d(x, y)$ for all $x, y \in M$. 
Thus, $F_{\gamma}$ is Lipschitz continuous with Lipschitz constant $1$ so that 
$F_{\gamma}$ is differentiable except for a measure zero set. 
Furthermore, on the differentiability of $F_{\gamma}$, 
we already know the following fact in general\,(cf.\,\cite[Theorem 1.1]{S2})\,: 
\begin{center}
$F_{\gamma}$ is differentiable at $x \in M$, 
if $x$ is {\em an interior point} of some ray $\sigma$ asymptotic to $\gamma$. 
\end{center} 
However, the differentiability at $\sigma(0)$ of $F_{\gamma}$ has never been studied in general. 

\bigskip

Our purpose in this section is to show the differentiability of a Busemann function 
of a ray at a starting point of some ray asymptotic to the ray (Theorem \ref{thm3.3}).

\bigskip

We will first recall the following fundamental property of a Busemann function\,:

\begin{lemma}{\rm (cf.\,\cite[Theorem 1.1]{S2})}\label{lemma3.2} 
For any ray $\sigma$ asymptotic to a ray $\gamma$ in $M$, 
\[
F_{\gamma} \circ \sigma (t) = t + F_{\gamma} \circ \sigma (0)
\]
holds for all $t \ge 0$.
\end{lemma}

By Lemma \ref{lemma3.2}, it is easy to prove that 

\begin{lemma}\label{lemma10091801}
If $F_{\gamma}$ is differentiable at a point $q \in M$, then
there exists a unique asymptotic ray $\sigma$ of $\gamma$ emanating from $q = \sigma (0)$. 
\end{lemma}

\begin{lemma}\label{lemma3.3}
If there exists a unique asymptotic ray $\sigma$ of $\gamma$ emanating from $q = \sigma (0)$, 
then $F_{\gamma}$ is differentiable at a point $q \in M$.
\end{lemma}

\begin{proof} 
To see this lemma, it is sufficient to prove the following relationship: 
\[
\lim_{s \to 0} \frac{F_{\gamma} (\tau(s)) - F_{\gamma} (q)}{s} = 
\left\langle \tau' (0), \sigma' (0) \right\rangle = 
\cos \left( \angle (\tau' (0), \sigma' (0)) \right)_.
\]
for all geodesic segments $\tau(s)$ emanating from $q$ 
in a convex neighborhood $U_{q}$ around $q$.\par 
Take a point $r \in U_{q}$ on the unique asymptotic ray $\sigma$ of $\gamma$. 
We define 
\[
\varphi (s) := d(r, \tau(s)) 
\]
for any fixed minimal geodesic segment $\tau (s)$ emanating from $q$ in $U_{q}$. 
It follows from the Taylor expansion of $\varphi (s)$ about $s = 0$ 
that there exists a constant $C_{1}$ such that 
\begin{align}\label{B1-0}
\varphi (s) 
& = \varphi (0) + \varphi'(0)s + \frac{1}{2!} \varphi''(0)s^{2} + \cdots \notag\\ 
& \le \varphi (0) + \varphi'(0)s + C_{1} {s^{2}}_.
\end{align}
By (\ref{B1-0}) and the first variation formula, we have  
\begin{align}\label{B1}
d(r, \tau(s)) = \varphi (s)
&\le \varphi (0) + \varphi'(0)s + C_{1} s^{2} \notag\\
&= d(r, \tau (0)) + s \cos  \left( \pi - \angle (\tau' (0), \sigma' (0)) \right) + C_{1}s^{2} \notag\\ 
&= d(r, q) - s \cos  \left(\angle (\tau' (0), \sigma' (0)) \right) + C_{1}{s^{2}}_.
\end{align}
Since $F_{\gamma} (r) - F_{\gamma} (q) = d(q, r)$ by Lemma \ref{lemma3.2}, 
we have, by (\ref{B1}),  
\begin{align}\label{B2}
F_{\gamma} (\tau(s)) - F_{\gamma} (q)
&= F_{\gamma} (\tau(s)) - F_{\gamma} (r)  + d(q, r)\notag\\ 
&\ge d(q, r) - d(r, \tau(s))\notag\\ 
&\ge s \cos  \left(\angle (\tau' (0), \sigma' (0)) \right) - {C_{1}s^2}_.
\end{align}
Thus, by (\ref{B2}), we have 
\begin{equation}\label{B3}
\liminf_{s \to 0} \frac{F_\gamma (\tau(s)) - F_\gamma (q)}{s} 
\ge  
\cos \left( \angle (\tau' (0), \sigma' (0) ) \right)_.
\end{equation}
On the other hand, for each sufficiently small $s$, 
let $\sigma_{s}$ be 
an asymptotic ray of $\gamma$ emanating from $\tau (s)$, 
and set $\eta_{s} := \angle (\sigma'_{s} (0), \tau' (s))$ and $\delta := d(q, r)$. 
Furthermore, we define, for fixed $s$, 
\[
\psi_{s}(t) : = d(\sigma_{s} (\delta), \tau(t))_.
\] 
Since the distance function $d(\, \cdot \, , \, \cdot \, )$ is $C^{\infty}$ around $(q, r)$, 
it follows from the Taylor expansion of $\psi_{s}(t)$ about $t = s$ 
that there exists a constant $C_{2}$ such that
\begin{align}\label{B4-0}
\psi_{s}(t)
& = \psi_{s}(s) + \psi'_{s}(s) (t - s) + \frac{1}{2!} \psi''_{s}(s) (t - s)^{2} + \cdots \notag\\ 
& \le \psi_{s}(s) + \psi'_{s}(s) (t - s) + C_{2} (t - s)^{2}. 
\end{align}
Substituting $t = 0$ for (\ref{B4-0}), 
we have, by the first variation formula,  
\begin{align}\label{B4}
d(\sigma_{s} (\delta), q) = d(\sigma_{s} (\delta), \tau(0)) = \psi_{s}(0) 
& \le  \psi_{s}(s) -s \psi'_{s}(s) + C_{2}s^{2} \notag\\ 
&\le d(\sigma_{s} (\delta), \tau(s)) -s \cos \left( \pi - \eta_{s} \right) + C_{2}s^2 \notag\\ 
&= d(\sigma_{s} (\delta), \sigma_{s} (0)) + s \cos  \eta_{s} + C_{2}s^2 
\end{align}
for each sufficiently small $s$. 
Since $\sigma_{s}$ is an asymptotic ray of $\gamma$ emanating from $\tau (s)$, 
we have, by Lemma \ref{lemma3.2}, 
\begin{equation}\label{B5-0}
F_{\gamma} (\sigma_{s} (\delta)) 
= \delta + F_{\gamma} (\sigma_{s} (0))
= \delta + F_{\gamma} (\tau(s)).
\end{equation}
By (\ref{B4}) and (\ref{B5-0}), we have  
\begin{align}\label{B5}
F_{\gamma} (\tau(s)) - F_\gamma (q)
&= F_{\gamma} (\sigma_{s} (\delta)) - F_{\gamma} (q)  -\delta\notag\\ 
&\le d(q, \sigma_{s} (\delta)) - d(q, \sigma(\delta))\notag\\ 
&= d(q, \sigma_{s} (\delta)) - d(\sigma_{s} (0), \sigma_{s} (\delta))\notag\\ 
&\le s \cos  \eta_{s} + C_{2}s^2
\end{align}
for each sufficiently small $s$. 
Thus, since $\sigma$ is the unique asymptotic ray of $\gamma$ emanating from $q = \sigma (0)$, 
we have, by (\ref{B5}),  
\begin{equation}\label{B6}
\limsup_{s \to 0} \frac{F_{\gamma} (\tau(s)) - F_{\gamma} (q)}{s} 
\le 
\lim_{s \to 0} \cos  \eta_{s} = \cos \left( \angle (\tau' (0), \sigma' (0) ) \right)_.
\end{equation}
By (\ref{B3}) and (\ref{B6}), we get the desired assertion. 
$\qedd$
\end{proof}

Furthermore, applying the method in \cite[Lemma 2.1]{IT} and \cite[Section 5]{MT}, 
we have the key lemma to prove Theorem \ref{thm3.3}.

\begin{lemma}\label{lem2.6}
Let $F : U \lra \R$ be a Lipschitz function defined on an open subset $U$ of 
Euclidean $n$-dimensional space $\R^{n}$. 
Then, a necessary and sufficient condition for $F$ to be differentiable at $x_{0} \in U$ is that 
the directional derivative function 
\[
F'_{+} (x_{0}  ; v) := \lim_{t \downarrow 0} \frac{F(x_{0} + t v) - F(x_{0})}{t}
\]
exists for all $v \in \R^{n}$ and is linear.  
\end{lemma}

\begin{proof}
It is sufficient to prove that 
$F$ is differentiable at $x_{0} \in U$, 
if the directional derivative function $F'_{+} (x_{0}  ; v)$ of $F$ exists for all $v \in \R^{n}$ and is linear. 
Under this assumption, we have 
\begin{equation}\label{L1}
\lim_{t \downarrow 0} {\frac{\left| F(x_{0} + t v) - F(x_{0}) - F'_{+} (x_{0}  ; tv) \right|}{t} = 0}
\end{equation}
for all unit vectors $v \in U$.  
Take any fixed $\ve > 0$. 
Since $\Sph^{n - 1}_{x_{0}} := \{v \in T_{x_{0}} (\R^{n}) \  | \ \| v \| = 1 \}$ is compact, 
we can take finitely many $v_{1},  \ldots , v_{k(\ve)} \in \Sph^{n - 1}_{x_{0}}$ 
such that for any $v \in \Sph^{n - 1}_{x_{0}}$, $\| v_{i} - v \| < \ve$ holds for some 
$i \in \{1, \ldots, k(\ve) \}$. 
Let $w$ be any nonzero vector at $x_{0}$, 
and $v_{i_{0}}$ be one of the vectors $v_{1}, \ldots , v_{k(\ve)}$ such that 
\[
{\left\| \left. \frac{1}{\| w \|} w - v_{i_{0}} \right. \right\| < \ve}_.
\] 
By the triangle inequality, we have 
\begin{align}\label{L2}
&\left| F (x_{0} + w) - F (x_{0}) - F'_{+} (x_{0}  ; w) \right| \notag\\
&\le \left| F (x_{0} + w) - F (x_{0} + \| w \| v_{i_{0}}) \right| 
+ \left| F (x_{0} + \| w \| v_{i_{0}}) - F (x_{0}) - F'_{+} (x_{0}  ; \| w \| v_{i_{0}}) \right| \\
&\mbox{} \quad 
+ \left| F'_{+} (x_{0}  ; w) - F'_{+} (x_{0}  ; \| w \| v_{i_{0}}) \right|_. \notag
\end{align}
Now, since $F$ is Lipschitz, 
there exists a Lipschitz constant number $C_{0}$ such that 
\begin{equation}\label{L3}
 \left| F (x_{0} + w) - F (x_{0} + \| w \| v_{i_{0}}) \right| 
 \le C_{0} \left\| \left. w - \| w \| v_{i_{0}} \right. \right\| < C_{0} \| w \| \ve_.  
\end{equation}
Furthermore, since $F'_{+} (x_{0}  ;  \, \cdot \, )$ is linear, 
\begin{align}\label{L4}
\left| F'_{+} (x_{0}  ; w) - F'_{+} (x_{0}  ; \| w \| v_{i_{0}}) \right| 
& \le C_{0} \left\| \left. w - \| w \| v_{i_{0}} \right. \right\| < C_{0} \| w \| \ve_.
\end{align}
By (\ref{L2}), (\ref{L3}), and (\ref{L4}), we see
\begin{align}\label{L5}
&\left| F (x_{0} + w) - F (x_{0}) - F'_{+} (x_{0}  ; w) \right| \notag\\
&< 2 C_{0} \| w \| \ve 
+ 
\left| F (x_{0} + \| w \| v_{i_{0}}) - F (x_{0}) - F'_{+} (x_{0}  ; \| w \| v_{i_{0}}) \right|_, 
\end{align}
for all nonzero vectors $w$ at $x_{0}$. 
Thus, by (\ref{L1}) and (\ref{L5}), we get 
\begin{equation}\label{L6}
\limsup_{\| w \| \to 0}
{ 
\frac{\left| F (x_{0} + w) - F (x_{0}) - F'_{+} (x_{0}  ; w) \right|}{\| w \|} 
\le 
2C_{0} \ve
}_.
\end{equation}
Since $\ve$ is arbitrarily taken, this inequality (\ref{L6}) shows that 
$F$ is differentiable at $x_{0}$. 
$\qedd$
\end{proof}

\begin{theorem}\label{thm3.3}
Let $F_{\gamma}$ be a Busemann function of a ray $\gamma$
on a complete non-compact Riemannian manifold $M$. 
Then, $F_{\gamma}$ is differentiable at a point $q \in M$ 
if and only if there exists a unique asymptotic ray $\sigma$ of 
$\gamma$ emanating from $q = \sigma (0)$. 
Moreover, the gradient vector of
$F_{\gamma}$ at a differentiable point $q$ equals the velocity vector of
the unique asymptotic ray of $\gamma$. 
\end{theorem}

\begin{proof}
The first assertion in this theorem is a direct consequence 
of Lemmas \ref{lemma10091801} and \ref{lemma3.3}. 
Finally, we will prove the second assertion in this theorem. 
It follows from Lemma \ref{lemma3.3} that 
\begin{equation}\label{B7}
dF_{\gamma} (v) = \lim_{s \to 0} \frac{F_{\gamma} \circ \exp_{q} (s v) - F_{\gamma} (q)}{s} = 
\left\langle v, \sigma' (0) \right\rangle
\end{equation}
for all $v \in T_{q}M$. 
Thus, by (\ref{B7}) and Lemma \ref{lem2.6}, we have 
\[
(\nabla F_\gamma)_{\sigma(0)} = \sigma'(0)
\]
for the unique asymptotic ray $\sigma$ of $\gamma$. 
$\qedd$
\end{proof}

\section{Proof of Main Theorem}\label{pmt}
\noindent

In this section, 
let $(M,p)$ denote a complete non-compact Riemannian $n$-manifold $M$
whose radial curvature at the base point $p$ is bounded from below by
that of a non-compact von Mangoldt surface of revolution $(\wt{M}, \tilde{p})$ 
with its total curvature $c(\wt{M}) > \pi$.\par
Let $A_{\tilde{q}}$ be the set of all unit vectors tangent to the rays emanating 
from a point $\tilde{q} \in \wt{M}$. 
Furthermore, we denote by $\mu$ the Lebesgue measure on the unit circle 
$\Sph^{1}_{\tilde{q}} := \{v \in T_{\tilde{q}} \wt{M} \ | \ \| v \| = 1 \}$ 
with total measure $2 \pi$.
Since $c(\wt{M}) > \pi$,  it follows from  \cite[Theorem 6.2.1]{SST} that 

\begin{lemma}\label{lem3.1}
There exist numbers $R, \delta > 0$ such that 
\[
\mu(A_{\tilde{q}})< \pi - 2\delta
\] 
for all $\tilde{q} \in \wt{M} \setminus B_{R}(\tilde{p})$. 
Here, 
$B_{R}(\tilde{p}) \subset \wt{M}$
is the open distance $R$-ball around $\tilde{p} \in \wt{M}$.
\end{lemma}

As a direct consequence of Lemma \ref{lem3.1}, we have the following lemma\,: 

\begin{lemma}\label{lemma4.3}
Let $R, \delta > 0$ be the numbers in Lemma \ref{lem3.1}. 
For each $\tilde{q} \in \wt{M} \setminus \overline{B_{R}(\tilde{p})}$, 
there exists a number $R_{1} > R$ such that 
\[
\angle (\tilde{p}\tilde{q}\tilde{x}) \ge \frac{\pi}{2} + \delta
\]
for all $\tilde{x} \in \wt{M} \setminus B_{R_{1}}(\tilde{p})$. 
Here, $\angle (\tilde{p}\tilde{q}\tilde{x})$ is the angle between minimal geodesics from $\tilde{q}$ to $\tilde{p}$ and $\tilde{x}$. 
\end{lemma}

\begin{remark}
Lemmas \ref{lem3.1} and \ref{lemma4.3} hold for all arbitrary 
non-compact model surfaces of revolution $\wt{M}$ admitting $c (\wt{M}) > \pi$ 
(See \cite[Section 2]{KT2}).
\end{remark}

In the proof of Lemma \ref{lem3.6}, we will apply a new type of the Toponogov comparison theorem. 
The comparison theorem was established by the present authors 
as a generalization of the comparison theorem in conventional comparison geometry, 
which is stated as follows\,:

\medskip

\begin{theorem}\label{TCT}\hspace{-1.5mm}{\rm (\cite[Theorem 4.12]{KT1})} 
Let $(M,p)$ be a complete non-compact Riemannian manifold $M$ 
whose radial curvature at the base point $p$ is bounded from below by
that of a non-compact model surface of revolution $(\wt{M}, \tilde{p})$. 
If $(\wt{M}, \tilde{p})$ admits a sector $\wt{V}(\delta_{0})$, 
$\delta_{0} \in (0, \pi]$, having no pair of cut points, 
then, for every geodesic triangle $\triangle(pxy)$ in $(M,p)$ 
with $\angle (xpy) < \delta_{0}$, 
there exists a geodesic triangle 
$\wt{\triangle} (pxy) :=\triangle(\tilde{p}\tilde{x}\tilde{y})$ 
in $\wt{V}(\delta_{0})$ such that
\[
d(\tilde{p},\tilde{x})=d(p,x), \quad d(\tilde{p},\tilde{y})=d(p,y), \quad d(\tilde{x},\tilde{y})=d(x,y) 
\]
and that
\[
\angle (xpy) \ge \angle (\tilde{x}\tilde{p}\tilde{y}), \quad  
\angle (pxy) \ge \angle (\tilde{p}\tilde{x}\tilde{y}), \quad
\angle (pyx) \ge \angle (\tilde{p}\tilde{y}\tilde{x}). 
\]
Here $\angle(pxy)$ denotes the angle between the minimal geodesic segments 
from $x$ to $p$ and $y$ forming the triangle $\triangle(pxy)$.
\end{theorem}

\begin{remark}
In \cite{KT4}, the present authors very recently generalized, from radial curvature 
geometry's standpoint, the Toponogov comparison theorem to a complete Riemannian manifold 
with smooth convex boundary.
\end{remark}

\begin{lemma}\label{lem3.6}
Let $R, \delta > 0$ be the numbers in Lemma \ref{lem3.1}. 
For any ray $\gamma$ on $M$ and any asymptotic ray
$\sigma$ of $\gamma$ emanating from a point 
$q = \sigma(0) \in M \setminus \overline{B_{R}(p)}$, we have 
\[
\angle(\sigma'(0), \alpha'(d(p, q))) \le \frac{\pi}{2} - \delta
\]
for all minimal geodesic segments $\alpha : [0, d(p, q)] \lra M$ emanating from $p$ to $q$. 
Here, $\overline{B_{R}(p)} \subset M$ is the closed distance $R$-ball around $p \in M$.
\end{lemma}

\begin{proof} 
For any sufficiently small fixed number $\ve > 0$, 
let $\alpha_{s}$ be a minimal geodesic emanating from 
$p$ to $\sigma(s)$ for each $s \in (0, \ve)$, and 
let $t_{0} > 0$ be a sufficiently large number such that
\[
d(p, \gamma(t_{0})) \gg \max_{s \in [0, \ve]} d(p, \sigma (s)).
\]  
Furthermore, for each $s \in (0, \ve)$, 
let $\sigma_{s, t_{0}}$ be a minimal geodesic emanating 
from $\sigma(s)$ to $\gamma(t_{0})$. 
Since $\wt{M}$ is a von Mangoldt surface of revolution, $\wt{M}$ has no pair of cut points in its 
sector $\wt{V} (\pi)$. 
Hence, we may apply Theorem \ref{TCT} to the geodesic triangle 
$\triangle (p \sigma (s) \gamma(t_{0}))$ 
consisting of $\alpha_{s}$, $\sigma_{s, t_{0}}$, and a minimal geodesic joining 
$p$ to $\gamma(t_{0})$. 
Thus, by the comparison theorem and Lemma \ref{lemma4.3}, 
\begin{equation}\label{eq:MT1.3}
\angle (- \alpha'_{s} (d(p, \sigma(s))), \sigma'_{s, t_{0}} (0)) \ge 
\angle (- \tilde{\alpha}'_{s} (d(p, \sigma(s))), \tilde{\sigma}'_{s, t_{0}} (0)) \ge \frac{\pi}{2} + \delta
\end{equation}
holds for each $s \in (0, \ve)$. 
Here, $ \tilde{\alpha}_{s}$ and $\tilde{\sigma}_{s, t_{0}}$ denote the edges of a geodesic triangle 
$\wt{\triangle} (p \sigma (s) \gamma(t_{0})) \subset \wt{M}$ 
corresponding to $\alpha_{s}$ and $\sigma_{s, t_{0}}$, respectively. 
On the other hand, since, for each $s \in (0, \ve)$, the sub\,-\,ray $\sigma|_{[s, \infty)}$ 
is the unique asymptotic ray of $\gamma$ emanating from $\sigma (s)$\,(cf.\,\cite[Theorem 1.1]{S2}), 
each $\sigma_{s, t_{0}}$ converges to the sub\,-\,ray $\sigma|_{[s, \infty)}$ 
as $t_{0} \to \infty$. 
Since $t_{0}$ is taken sufficiently large, we see, by (\ref{eq:MT1.3}), for each fixed $s \in (0, \ve)$,   
\begin{equation}\label{eq:MT1.3.1}
\angle (- \alpha'_{s} (d(p, \sigma(s))), \sigma'(s)) \ge {\frac{\pi}{2} + \delta}_.
\end{equation}
Moreover, since $\Sph^{n - 1}_{p} := \left\{ v \in T_{p} M \ | \ \| v \| = 1 \right\}$ is compact, 
there exists a sequence $\{ s_{n} \} \to 0$ as $n \to \infty$ 
such that $\alpha_{s_n}$ converges to some minimal geodesic 
$\alpha_{0}$ emanating $p$ to $q$. 
Then, since $\sigma' (s_{n})$ converges to $\sigma'(0)$ 
as $n \to \infty$, we have, by (\ref{eq:MT1.3.1}), 
\begin{equation}\label{eq:MT1.4}
\angle (- \alpha'_{0} (d(p, q)), \sigma' (0)) \ge {\frac{\pi}{2} + \delta}_.
\end{equation}
Thus, by (\ref{eq:MT1.4}) and \cite[Lemma 2.1]{IT}, 
we have, for any minimal geodesic segment $\alpha : [0, d(p, q)] \lra M$ emanating from $p$ to $q$,  
\begin{equation}\label{eq:MT1.5}
\angle (- \alpha' (d(p, q)), \sigma' (0)) \ge 
\angle (- \alpha'_{0} (d(p, q)), \sigma' (0)) \ge {\frac{\pi}{2} + \delta}_.
\end{equation}
Therefore, we get, by this (\ref{eq:MT1.5}), 
\[
\angle(\sigma'(0), \alpha'(d(p, q))) = \pi - \angle (- \alpha' (d(p, q)), \sigma' (0)) \le {\frac{\pi}{2} - \delta}_,
\]
which is the desired assertion in this lemma. 
$\qedd$
\end{proof}

Recall that a point $x \in M$ is called a {\it critical point for} 
a Busemann function $F_{\gamma}$ of a ray $\gamma$ on $M$ 
if, for every nonzero unit vector $v \in T_{x} M$, 
there exists an asymptotic ray $\sigma$ of $\gamma$ emanating from $\sigma(0) = x$
such that $\angle(v, \sigma'(0)) \le \pi/2$\,(cf.\,\cite{S3}). 
Therefore, as a direct consequence of Lemma \ref{lem3.6}, 
we have the following theorem, which shows the (MT--2) in the Main Theorem\,:  

\begin{theorem}\label{prop3.6-1}
$\overline{B_{R}(p)}$ contains all critical points for any Busemann function on $M$.    
\end{theorem}

Hereafter, we are going to prove Main Theorem.  

\begin{lemma}\label{lem2.7}
Let $F : [ a, b ] \lra \R$ be a Lipschitz function. 
Then its derivative function $F'(x)$ is integrable and
\begin{equation}\label{L7}
F(x)-F(a)=\int_a^x F ' (t) dt
\end{equation}
holds for all $x \in [ a, b ]$.
\end{lemma}

\noindent
The proof of Lemma \ref{lem2.7} is found, for example, in \cite[Theorem 7.29]{WZ}.

\begin{remark}
Lipschitz continuity is a very important property. For example, the
Cantor--Lebesgue function admits an integrable derivative function, i.e., $0$, 
but this function does not satisfy the equation (\ref{L7}).
\end{remark}

\begin{lemma}\label{lem3.7}
Let $F_{\gamma}$ be a Busemann function of a ray $\gamma$ on $M$ 
and $\alpha(t)$ a minimal geodesic segment in $M$ emanating from $p$. 
If $F_\gamma$ is differentiable at $\alpha(t)$ for almost all $t \in (a, b)$, 
and satisfies 
\[
{\angle ((\nabla F_\gamma)_{\alpha(t)}, \alpha'(t))\leq \frac{\pi}{2}-\delta}_, 
\] 
then, we have  
\[
F_\gamma(\alpha(b))-F_\gamma(\alpha(a)) \ge (b - a) \sin\delta.
\]
\end{lemma}

\begin{proof}
By assumptions on this lemma and Theorem \ref{thm3.3}, we have 
\begin{align}\label{eq:MT2.1}
\frac{d}{dt}F_{\gamma}(\alpha(t)) 
&= \langle (\nabla F_\gamma)_{\alpha(t)}, \alpha'(t) \rangle \notag\\ 
&= \cos \left( \angle ((\nabla F_\gamma)_{\alpha(t)}, \alpha'(t)) \right) \notag\\
&\ge \cos \left( \frac{\pi}{2}-\delta \right) \notag\\
&= \sin \delta 
\end{align}
for almost all $t\in(a, b)$.
Therefore, by Lemma \ref{lem2.7} and (\ref{eq:MT2.1}), we get 
 \begin{align*}
F_\gamma(\alpha(b))-F_\gamma(\alpha(a)) 
= \int^{b}_{a} \frac{d}{dt}F_{\gamma}(\alpha(t)) dt  
\ge (b - a) \sin \delta. 
\end{align*}
$\qedd$
\end{proof}

\begin{lemma}\label{lem3.4}
Let $E \subset \R^n$ be a subset of Lebesgue measure zero.
Then for almost all $x \in \R^{n-1},$
the set 
\[
E_{x} : = \{ t \in \R \ | \ (t, x) \in E \}
\]
is a subset of $\R$ of Lebesgue measure zero.
\end{lemma}

\noindent
The proof of Lemma \ref{lem3.4} is written in \cite[Lemma 6.5]{WZ}.

\begin{lemma}\label{lem3.5} {\bf(Rademacher)}
If a function $F : U \lra \R$ is Lipschitz on an open subset $U$ of $\R^n$, 
then $F$ is differentiable almost everywhere.
\end{lemma}
 
\noindent
One may find the proof of this Rademacher theorem in \cite{Mo}. 

\begin{lemma}\label{pro3.8}
For any Busemann function $F_\gamma$ on $M$ and any 
$q \in M \setminus \overline{B_{R}(p)}$, we have  
\[
F_\gamma(q)-F_\gamma(\alpha(R)) \ge \left( d(p,q) - R \right) \sin \delta.
\]
Here $\alpha : [ 0, d(p,q) ] \lra M$ denotes a minimal geodesic segment joining $p$ to $q$.
\end{lemma}

\begin{proof}
It is sufficient to prove that, for each $t_0 \in(R, d(p, q))$, 
there exists a number $\ve_0 > 0$ such that
\begin{equation}\label{20100919eq1}
F_{\gamma}(\alpha(t))-F_{\gamma}(\alpha(s)) \geq (t - s ) \sin \delta
\end{equation}
for all $t_0 - \ve_0 < s < t < t_0 + \ve_0$. 
Since $\alpha(t_0)$ is not a cut point of $p = \alpha(0)$,  
there exist an open neighborhood $\cU$ around 
$\alpha' (0)$ in $\Sph^{n - 1}_{p}$, 
an open neighborhood $V$ around $\alpha(t_0)$ in $M$ 
and an open interval $(t_0 - \ve_0, t_0 + \ve_0)$ such that
$\cU \times (t_0 - \ve_0, t_0 + \ve_0)$ is diffeomorphic to $V$ 
by a map $\varphi$, where $\varphi^{-1}(v, t) := \exp_{p} (t v)$.  
Since $F_\gamma \circ \varphi^{-1}$
is Lipschitz, it follows from Lemma \ref{lem3.5} that 
there exists a subset $\cE$ of $\R^n$ of Lebesgue measure zero such that 
$F_\gamma \circ \varphi^{-1}$ is differentiable on 
$(t_0 - \ve_0, t_0 + \ve_0) \times \cU \setminus \cE$. 
Here, $\cU$ is identified as an open subset of $\R^{n-1}$. 
By Lemma \ref{lem3.4}, 
there exists a sequence $\{ \alpha_{i} \}$ of minimal geodesic segments emanating
from $\alpha(0)$ converging to $\alpha$ such that $F_{\gamma}$ is differentiable 
almost everywhere on $\alpha_{i} (t_{0} - \ve_0, t_{0} + \ve_0)$ for each $i.$
By Lemmas \ref{lem3.6} and \ref{lem3.7}, 
we have 
$F_{\gamma}(\alpha_{i} (t) ) - F_{\gamma}(\alpha_{i} (s) ) \ge ( t - s ) \sin \delta$ 
for all $t_{0} - \ve_0 < s < t < t_{0} + \ve_0$ and all $i$. 
By taking the limit, we get (\ref{20100919eq1}).
$\qedd$
\end{proof}

\begin{theorem}\label{20100919thm1}
All Busemann functions $F_{\gamma}$ on $M$ are exhaustions.  
\end{theorem}

\begin{proof}
Set 
$M(F_{\gamma} : a) := \left\{ x \in M  \ | \ F_{\gamma} (x) \le a \right\}$
for each $a \in \R$. 
Since $\partial B_{R}(p)$ is compact in Lemma \ref{pro3.8}, 
there exists a constant $\cN(R)$ such that 
$F_{\gamma} (\alpha (R)) \ge \cN(R)$. 
By Lemma \ref{pro3.8}, we have 
\begin{equation}\label{F4}
a - \cN(R) \ge F_{\gamma} (q) - F_{\gamma} (\alpha (R)) \ge \left( d(p, q) - R \right) \sin \delta
\end{equation}
for all $q \in ( M \setminus \overline{B_{R}(p)} ) \cap M(F_{\gamma} : a)$. 
By (\ref{F4}), we get,  
\begin{equation}\label{F5}
d(p, q) \le {\frac{a - \cN(R)}{\sin \delta} + R}
\end{equation}
for all $q \in ( M \setminus \overline{B_{R}(p)} ) \cap M(F_{\gamma} : a)$. 
Thus, (\ref{F5}) implies that $M(F_{\gamma} : a)$ is compact. 
$\qedd$
\end{proof}

 \section{Proof of Corollary}

\begin{corollary} Let $(M,p)$ be a complete non-compact Riemannian $n$-manifold $M$ 
whose radial curvature at the base point $p$ is bounded from below by
that of a non-compact von Mangoldt surface of revolution $(\wt{M}, \tilde{p})$.
If $(\wt{M}, \tilde{p})$ admits $c (\wt{M}) > \pi$,
\begin{enumerate}[{\rm ({C\,--\,}1)}]
\item
$M$ has finite topological type, that is, 
$M$ is homeomorphic to the interior of a compact manifold with boundary. 
\item
The isometry group of $M$ is compact. 
\end{enumerate}
\end{corollary}

\begin{proof}
We first prove (C\,--\,1). 
Recall that, for a fixed point $q \in M$, 
a point $x \in M \setminus \{ q \}$ is called a {\it critical point of} $d(q, \, \cdot \, )$ 
(or {\em critical point for} $q$) if, 
for every nonzero tangent vector $v \in T_{x} M$, 
we find a minimal geodesic $\gamma$ emanating from $x$ to $q$
satisfying $\angle(v, \gamma'(0)) \le \pi/2$ (See \cite{GS}). 
Thus, Lemma \ref{lem3.6} implies that there exists a number $\ve > 0$ such that 
$M \setminus B_{R + \ve}(p)$ has no critical point of  $d(p, \, \cdot \, )$, 
where $R > 0$ is the number in Lemma \ref{lem3.1}. 
By \cite[Isotopy Lemma]{Gv}, $M \setminus B_{R + \ve}(p)$ is homeomorphic to 
$\partial B_{R + \ve}(p) \times [R + \ve, \infty)$ so that $M$ has finite topological type.\par 
Next, we prove (C\,--\,2). 
Let $I(M)$ be the isometry group of $M$. 
We will first recall Myers and Steenrod' method in \cite[Section 7]{MS}. 
They took $n + 1$ points $x_0, x_1, \ldots, x_n$ in $M$ 
which are independent in a certain sense 
(roughly speaking, $x_1, x_2, \ldots, x_n$ are different points to each other 
in a convex neighborhood around $x_0$ with 
$x_i \not= x_0$ for $i = 1, 2, \ldots, n$ such that each vector $v_{i}$ at $x_0$, which is  
tangent to each minimal geodesic emanating from $x_0$ to $x_i$ for $i = 1, 2, \ldots, n$, 
is linearly independent), 
and proved that the mapping 
\[
g \in I(M) \lra \left( g(x_0), g(x_1), \ldots, g(x_n) \right) \in M^{n + 1}
\]
is one\,-to\,-\,one and has a closed submanifold of $M^{n + 1}$ as its image. 
Here, 
$M^{n + 1}$ is the product manifold $M \times M \times \cdots \times M$ 
taken with itself $n + 1$ times. 
Thus, it is sufficient to prove that $d(p , g(q))$ is bounded from above by some constant 
for all $g \in I(M)$ and for some point $q \in M$ with a certain special property. 
Now, let $F_{\gamma}$ be a Busemann function of any ray $\gamma$ on $M$. 
Since $F_{\gamma}$ is an exhaustion by Theorem \ref{20100919thm1}, 
$F_{\gamma} ( \, \cdot \, )$ attains its minimum at some $q_{0} \in M$. 
Thus, $q_{0}$ is a critical point for $F_{\gamma}$, 
and hence $q_{0} \in \overline{B_{R}(p)}$ by Theorem \ref{prop3.6-1}, 
where $R > 0$ is the number in Lemma \ref{lem3.1}. 
For any $g \in I(M)$, we see $g(q_0) \in \overline{B_{R}(p)}$, 
since $g(q_{0})$ is a critical point for $F_{g \circ \gamma}$ of the ray $g \circ \gamma$ on $M$. 
Thus, we get $d(p , g(q_{0})) \le R$ for all $g \in I(M)$. 
Therefore, $I(M)$ is compact.
$\qedd$
\end{proof}

\medskip

\bigskip

\begin{center}
Kei KONDO $\cdot$ Minoru TANAKA 

\bigskip
Department of Mathematics\\
Tokai University\\
Hiratsuka City, Kanagawa Pref.\\ 
259\,--\,1292 Japan

\bigskip
$\bullet$\,our e-mail addresses\,$\bullet$

\bigskip 
\textit{e-mail of Kondo}\,: 

\medskip
keikondo@keyaki.cc.u-tokai.ac.jp

\medskip
\textit{e-mail of Tanaka}\,: 

\medskip
m-tanaka@sm.u-tokai.ac.jp
\end{center}

\end{document}